\documentclass{elsart}

\usepackage{amssymb}
\usepackage{amsmath}

\begin{document}
\newcommand{\eop}{{\unskip\nobreak\hfil\penalty50
          \hskip2em\hbox{}\nobreak\hfil$\Box$
          \parfillskip=0pt \finalhyphendemerits=0 \par}}

\newtheorem{Theorem} {Theorem} [section]
\newtheorem{Proposition} [Theorem] {Proposition}
\newtheorem{Lemma} [Theorem] {Lemma}
\newtheorem{Corollary} [Theorem] {Corollary}
\newtheorem{Conjecture}[Theorem]{Conjecture}
\newtheorem{Example}[Theorem]{Example}

\numberwithin{equation}{section}

\begin{frontmatter}

\title{Semi-regular Relative Difference Sets with Large Forbidden Subgroups}
\author{Tao Feng\thanksref{supp1}}
\ead{ift@pku.edu.cn}
\address{School of Mathematical Sciences\\ Peking University, Beijing 100871, China}
\thanks[supp1]{Tao Feng was supported in part by National Natural Science Foundation of
China grant 10331030. He would like to thank Professor Weisheng Qiu
for supervision and encouragement.}
\author{Qing Xiang\thanksref{supp2}}\ead{xiang@math.udel.edu}
\address{Department of Mathematical Sciences\\ University of Delaware, Newark, DE 19716, USA}
\thanks[supp2]{Qing Xiang was supported in part by NSF grant DMS 0701049.}


\begin{abstract}
Motivated by a connection between semi-regular relative difference
sets and mutually unbiased bases, we study relative difference sets
with parameters $(m,n,m,m/n)$ in groups of non-prime-power orders.
Let $p$ be an odd prime. We prove that there does not exist a
$(2p,p,2p,2)$ relative difference set in any group of order $2p^2$,
and an abelian $(4p,p,4p,4)$ relative difference set can only exist
in the group $\Bbb{Z}_2^2\times \Bbb{Z}_3^2$. On the other hand, we
construct a family of non-abelian relative difference sets with
parameters $(4q,q,4q,4)$, where $q$ is an odd prime power greater
than $9$ and $q\equiv 1$ (mod $4$). When $q=p$ is a prime, $p>9$,
and $p\equiv$ 1 (mod 4), the $(4p,p,4p,4)$ non-abelian relative
difference sets constructed here are genuinely non-abelian in the
sense that there does not exist an abelian relative difference set
with the same parameters.
\end{abstract}

\end{frontmatter}

\section{Introduction}
Let $G$ be a finite (multiplicative) group of order $mn$, and let
$N$ be a subgroup of $G$ of order $n$. A $k$-subset $R$ of $G$ is
called an $(m,n,k,\lambda)$ {\it relative difference set} (RDS) in
$G$ relative to $N$ if every element $g\in G\setminus N$ has exactly
$\lambda$ representations $g=r_1 r_2^{-1}$ with $r_1$, $r_2\in R$,
and no non-identity element of $N$ has such a representation. The
subgroup $N$ is usually called {\it the forbidden subgroup}. If the
group $G$ is abelian (resp. non-abelian), then $D$ is called an {\it
abelian (resp. non-abelian) relative difference set}. When $n=1$,
$R$ is an $(m,k,\lambda)$ difference set in the usual sense. If
$k=n\lambda$, then $R$ is said to be {\it semi-regular}.

For a subset $X$ of $G$, we set $X^{(-1)}=\{x^{-1}\,|\, x\in X\}$;
also we use the same $X$ to denote the group ring element
$\sum_{x\in X}x\in \Bbb{Z}[G]$.  Then, a $k$-subset $R$ of $G$ is an
$(m,n,k,\lambda)$ relative difference set in $G$ relative to $N$ if
and only if it satisfies the following equation in the group ring
$\Bbb{Z}[G]$:
$$RR^{(-1)}=k+\lambda(G-N).$$

Character theory is a very useful tool in the study of difference
sets and relative difference sets in abelian groups. We state the
Fourier inversion formula below, which will be used many times in
the paper.

 \textbf{Inversion formula}\hspace{0.06in} Let $G$ be an abelian group of order $v$.
 If $A=\sum_{g\in G}a_g g\in \Bbb{Z}[G]$, then
 $a_h=\frac{1}{v}\sum_{\chi\in\hat G}\chi(Ah^{-1})$, for all $h\in G$, where
 $\hat{G}$ is the group of characters of $G$ and $\chi(Ah^{-1})=\sum_{g\in G}a_g\chi(gh^{-1})$.

One consequence of the inversion formula is as follows. Let $G$ be
an abelian group of finite order, and let $A$ and $B$ be two
elements of $\Bbb{Z}[G]$. Then we have $A=B$ if and only if
$\chi(A)=\chi(B)$ for all characters $\chi$ of $G$. The following
result is a standard characterization of relative difference sets by
using their character values (c.f. \cite[p.~374]{bjl}).

\begin{Proposition}\label{char}
Let $G$ be an abelian group of order $mn$ with a subgroup $N$ of
order $n$. Let $k$ and $\lambda$ be positive integers satisfying
$k(k-1)=\lambda n(m-1)$. Then a $k$-subset $D$ of $G$ is an
$(m,n,k,\lambda)$ difference set in $G$ relative to $N$ if and only
if for every non-principal character $\chi$ of $G$,
\begin{equation}
\chi(D)\overline{\chi(D)}= \left\{ \begin{array}{lll} k, & \mbox{if}
& \chi |_N \ne 1 \,
,\\
k-\lambda n, & \mbox{if} & \chi |_N = 1 \, .  \end{array} \right.
\end{equation}
where $\chi |_N$ is the restriction of $\chi$ to $N$.
\end{Proposition}

Recently a connection between semi-regular abelian RDS and mutually
unbiased bases is established in \cite{gr}. To explain the
connection, we first give the definition of mutually unbiased bases.
Let $\Bbb{C}$ be the field of complex numbers.  A pair of bases
$x_1,x_2,\ldots ,x_d$ and $y_1,y_2,\ldots ,y_d$ of $\Bbb{C}^d$ is
said to be {\it mutually unbiased} if they are both orthonormal and
there is a constant $a$ such that
$$|\langle x_i,y_j\rangle|=a,$$
for all $i$ and $j$, $i\neq j$, where $\langle\cdot,\cdot\rangle$ is
the standard inner product of $\Bbb{C}^d$. Let $N_{\rm MUB}(d)$
denote the maximum size of any set containing pairwise mutually
unbiased bases (MUB) of $\Bbb{C}^d$. It is an open question to
determine $N_{\rm MUB}(d)$ for every $d$. There are some
similarities between $N_{\rm MUB}(d)$ and $N_{\rm MOLS}(d)$, the
maximum number of mutually orthogonal Latin squares of size $d$. For
example, it is known \cite{dgs} that $N_{\rm MUB}(d)\leq d+1$; and
when $d=p^e$ is a prime power it was shown \cite{iv}, \cite{wf} that
$N_{\rm MUB}(p^e)=p^e+1$. Also if $d=st$, then we have
\begin{equation}\label{product}
N_{\rm MUB}(d)\geq {\rm min}\{N_{\rm MUB}(s), N_{\rm MUB}(t)\}.
\end{equation}
For an arbitrary positive integer $d$ and a prime $p$, we use
$\nu_p(d)$ to denote $p^{\alpha}$, where $p^{\alpha}|d$ but
$p^{\alpha+1}\nmid d$. We also use $\pi(d)$ to denote the set of
prime divisors of $d$. Then by (\ref{product}), we have
\begin{equation}\label{macNeish}
N_{\rm MUB}(d)\geq {\rm min}_{p\in \pi(d)}\{N_{\rm
MUB}(\nu_p(d))\}={\rm min}_{p\in \pi(d)}\{\nu_p(d) +1\}.
\end{equation}
We will refer to this construction as the {\it reduce to prime power
construction}. For more information on $N_{\rm MUB}(d)$, we refer
the reader to \cite{acw} and \cite{gr}.

We now state a theorem in \cite{gr} which establishes a connection
between semi-regular abelian RDS and mutually unbiased bases.

\begin{Theorem}\label{semiRDSMUB}{\rm (\cite{gr})}
The existence of a semi-regular $(m,n,m,m/n)$ RDS in an abelian group implies the existence of a set of $n+1$ mutually unbiased bases of $\Bbb{C}^m$.
\end{Theorem}

The proof of Theorem~\ref{semiRDSMUB} is a straightforward
application of Proposition\ref{char}. We refer the reader to
\cite{gr} for the proof and other background materials on mutually
unbiased bases. Motivated by the desire to use
Theorem~\ref{semiRDSMUB} to construct more MUB than the minimum in
(\ref{macNeish}) given by the reduce to prime power construction,
Wocjan \cite{pw} asked the following question: Does there exist an
abelian semi-regular relative difference set with parameters
$(m,n,m,m/n)$ satisfying
\begin{equation}\label{REQinq}
n> {\rm min}_{p\in \pi(m)}\{\nu_p(m)\}?
\end{equation}

We make some preliminary observations regarding this question. First
of all, most known semi-regular RDS have parameters
$(p^a,p^b,p^a,p^{a-b})$, where $p$ is a prime. The parameters of
these RDS will not satisfy (\ref{REQinq}). The reason is quite
straightforward. Note that if $m$ is a prime power, then
(\ref{REQinq}) simply becomes $n>m$. For RDS with parameters
$(p^a,p^b,p^a,p^{a-b})$, where $p$ is a prime, we have $p^{a-b}\geq
1$; hence $n=p^b\leq p^a=m$. Therefore to answer the question of
Wocjan we have to consider semi-regular $(m,n,m,m/n)$ RDS with $m$
not a prime power. As far as we know, there are only two general
constructions (\cite{mll}, \cite{djm}) of such semi-regular RDS with
$n>2$. The RDS constructed in these papers have parameters
\begin{equation}\label{nonprimepower}
(p^{2t}(p+1),p+1,p^{2t}(p+1),p^{2t}),
\end{equation}
where $t$ is a positive integer, and $p=2$ or $p$ is a Mersenne
prime. Note that the parameters in (\ref{nonprimepower}) do not
satisfy (\ref{REQinq}) either since $n=p+1$ and ${\rm
min}_{r\in\pi(m)}\{\nu_r(m)\}=p+1$ (here $p=2$ or $p$ is a Mersenne
prime).

Therefore we are motivated to search for semi-regular RDS with
parameters $(m,n,m,m/n)$ not of the form (\ref{nonprimepower}) and
$m$ not a prime power. The simplest case to consider is when
$(m,n,m,m/n)=(2p,p,2p,2)$, $p$ an odd prime. We prove in Section 3
that there does not exist a $(2p,p,2p,2)$ RDS in any group of order
$2p^2$. Next we prove that an abelian $(4p,p,4p,4)$ RDS with $p$ an
odd prime can only exist in the group $\Bbb{Z}_2^2\times
\Bbb{Z}_3^2$. On the construction side, we construct a family of
$(4q,q,4q,4)$ non-abelian RDS, where $q$ is an odd prime power
greater than $9$, $q\equiv 1$ (mod $4$). When $q=p$ is a prime (also
$p>9$ and $p\equiv 1$ (mod 4)), by the above nonexistence result on
abelian $(4p,p,4p,4)$ RDS, we see that the RDS we construct here are
genuinely nonabelian in the sense that there does not exist an
abelian RDS with the same parameters.

We give some preparation results in the rest of this section. For
any group $G$ with a subgroup $N$, we use $C_G(N)$ to denote the
{\it centralizer} of $N$ in $G$, namely, $C_G(N)=\{x\in G :
xy=yx,\;\forall y\in N\}$. Also we use exp$(G)$ to denote the
exponent of $G$. The following lemma on RDS is implicitly contained
in \cite{gal_ga}, and has its origin in \cite{pott_adds}.

\begin{Lemma}\label{lem_gp}
Let $G$ be a group of order $mn$ with an abelian normal subgroup $N$
of order $n$, and let $R$ be an $(m,n,m,m/n)$ RDS in $G$ relative to
$N$.  Then exp$(C_G(N))$ divides $2m$. Furthermore if the Sylow
$2$-subgroup of $N$ is not cyclic or $m/n$ is even, then
exp$(C_G(N))$ divides $m$.\end{Lemma}

\begin{pf} Since $N$ is abelian, we have $C_G(N)\geq N$. If
$C_G(N)=N$, then of course $|C_G(N)|=|N|$. Hence exp$(C_G(N))$
divides $|N|=n$, which in turn divides $m$ since $m/n$ is an
integer. So we will assume that $C_G(N)\neq N$ from now on. Given an
element $g\in G$, we use $\bar g$ to denote its image in $G/N$. Also
we use $r_{\bar g}$ to denote the unique element in $R\cap gN$. Now
for any given $g\in C_G(N)\setminus N$, we set
$$S=\{(r_{\overline{gh}}, r_{\bar h}) : \bar{h}\in G/N\}.$$
We have $|S|=m$. Since $N$ is normal in $G$, we see that for any
pair $(r_1,r_2)\in S$, $r_1r_2^{-1}\in gN$. Next we claim that each
$gu$, where $u\in N$, can be represented as $gu=r_1r_2^{-1}$, for
$m/n$ pairs $(r_1,r_2)\in S$. This claim can be seen as follows.
Since $R$ is an $(m,n,m,m/n)$ RDS in $G$ relative to $N$, each $gu$,
$u\in N$, can be represented as $gu=xy^{-1}$, for $m/n$ pairs
$(x,y)\in R\times R$. Let $y=hu'$, where $u'\in N$. Then
$x=guhu'=gh(h^{-1}uh)u'$. Since $N$ is normal in $G$, we have
$h^{-1}uh\in N$. Hence $x\in R\cap ghN$. The claim is proved. It
follows that,
$$g^m(\prod_{u\in N}u)^{m/n}=\prod_{u\in
N}(gu)^{m/n}=\prod_{(r_1,r_2)\in S}r_1r_2^{-1}.$$ Now using the
assumption that $g\in C_G(N)$, we can arrange the terms in the last
product above in such a way that $r_1r_2^{-1}$ is followed by
$r_2r_3^{-1}$, and so on. Therefore we have
$$g^m(\prod_{u\in N}u)^{m/n}=1.$$
The element $a:=\prod_{u\in N}u$ has order at most $2$. So
$g^{2m}=1$. Hence exp$(C_G(N))$ divides $2m$. If the Sylow
$2$-subgroup of $N$ is not cyclic, then $N$ has at least two
elements of order $2$; hence $a=1$. Therefore we have $g^m=1$ and
exp$(C_G(N))$ divides $m$. If $m/n$ is even, then clearly we have
$g^m=1$ and exp$(C_G(N))|m$. The proof is complete. \eop
\end{pf}

Let $p$ be a prime and $f:\Bbb{Z}_p^n\rightarrow \Bbb{Z}_p$ be a
function. The {\it Fourier transform} $\hat{f}$ of $f$ is defined by
$$\hat{f}({\mathbf b})=\sum_{{\mathbf x}\in\Bbb{Z}_p^n}\xi_p^{f({\mathbf x})+{\mathbf b}\cdot {\mathbf x}},\;\forall{\mathbf b}\in\Bbb{Z}_p^n,$$
where ${\mathbf b}\cdot {\mathbf x}$ is the standard dot product and
$\xi_p$ is a primitive $p$th root of unity in $\Bbb{C}$. The
function $f$ is said to be {\it $p$-ary bent} if $|\hat{f}({\mathbf
b})|=p^{n/2}$ for all ${\mathbf b}\in\Bbb{Z}_p^n$. In Section 4, we
will need the following theorem from \cite{hxd}.

\begin{Theorem}\label{hou} {\rm (\cite{hxd})}
Let $p$ be an odd prime. Then a function
$f:\Bbb{Z}_p\rightarrow\Bbb{Z}_p$ is $p$-ary bent if and only if
deg$(f)=2$.
\end{Theorem}

Throughout this paper, we fix the following notation: For a
multiplicative group $G$, we denote its identity by $1_G$, or simply
by $1$ if there is no confusion. For a positive integer $m$, $\xi_m
$ denotes a primitive $m$th root of unity in $\Bbb{C}$. For an odd
prime $p$, $\left(\frac{\cdot}{p}\right)$ is the Legendre symbol;
also we let
$$\Delta=\sum_{x\in \Bbb{Z}_p}\xi_p^{x^2}=\sum_{i=0}^{p-1}\left(\frac{i}{p}\right)\xi_p^i.$$
It is well known \cite{ir} that $\Delta\bar{\Delta}=p$ and
$\Delta=\pm\sqrt{p^*}$, where $p^*=(-1)^{\frac{p-1}{2}}p$. For an
integer $t$ such that $p\nmid t$, we use $\sigma_t$ to denote the
element in $Gal(\Bbb{Q}(\xi_p)/\Bbb{Q})$ that maps $\xi_p$ to
$\xi_p^t$. We have
$\sigma_t(\Delta)=\left(\frac{t}{p}\right)\Delta$. We will use
standard facts on prime ideal decompositions of rational integers in
cyclotomic fields freely. The readers are referred to \cite{was},
\cite{ir}, \cite{Mapl} for proofs of these facts.

\section{A construction of $(4q,q,4q,4)$ RDS in non-abelian groups}

In this section, we construct a family of $(4q,q,4q,4)$ RDS in
certain non-abelian groups of order $4q^2$, where $q$ is an odd
prime power, $q\equiv 1$ (mod $4$), and $q>9$.

For prime power $q=p^n$, $n\geq 1$, $p$ an odd prime, let
$K:=\Bbb{F}_q$ be the finite field of order $q$,
$K^*=K\setminus\{0\}$, and $tr: K\rightarrow \Bbb{F}_p$ be the
absolute trace function.  The quadratic character $\eta$ on $K$ is
defined by
\begin{align*}
\eta(x)=\begin{cases}1, & \text{if}\; x\; \text{is a nonzero square of}\; K,\\
                      0, & \text{if}\; x=0,\\
                      -1, & \text{if}\; x\; \text{is a nonsquare of}\; K.\end{cases}
\end{align*}
For $u\in K^*$, we define
$$S(u):=\sum_{x\in K}\xi_p^{tr(ux^2)}.$$
For simplicity, we write $S$ for $S(1)$. We have $S+S(u)=2\sum_{x\in
K}\xi_p^{tr(x)}=0$ if $u$ is a nonsquare of $K$. Therefore
$S(u)=\eta(u)S$ for every $u\in K^*$.

The {\it quadratic Gauss sum} $g(\eta)$ is defined by
$$g(\eta):=\sum_{x\in K}\eta(x)\xi_p^{tr(x)}.$$
Straightforward computations show that $g(\eta)=S$. Therefore
$$S\overline{S}=g(\eta)\overline{g(\eta)}=q,$$
c.f.~\cite[p.~11]{gjs}.

In the rest of this section we assume that $q\equiv 1$ (mod $4$),
$e$, $f$ are elements of $K$ satisfying $e^4=1$, $f^2=-1$,
respectively.

Given an element $s_2\in K^*$, we define
$$s_1=\frac{1}{2}((1+s_2)+\frac{f}{e^2}(1-s_2)),$$
$$s_3=\frac{1}{2}((1+s_2)-\frac{f}{e^2}(1-s_2)).$$

\begin{Lemma}\label{ss}
If $q>9$, then there exists $s_2\in K^*$ such that
$$\eta(s_1s_2s_3)=-1.$$
\end{Lemma}
\begin{pf} First, note that if $s_2\neq \frac{f+1}{f-1}$ or
$\frac{f-1}{f+1}$, then $s_1\neq 0$ and $s_3\neq 0$. Secondly,
$$s_1s_2s_3=\frac{s_2}{4}\left((1+s_2)^2-\frac{f^2}{e^4}(1-s_2)^2\right)=\frac{s_2}{2}(1+s_2^2).$$
Hence the number of $s_2\in K^*$ satisfying
$\eta(s_1s_2s_3)=\eta(2s_2(1+s_2^2))=-1$ is at least
\begin{equation}\label{quantity}\sum_{x\in K^*}\frac{1-\eta(2x(1+x^2))}{2}
-2=\frac{1}{2}\left(q-5-\sum_{x\in
K^*}\eta(2x+x^3)\right).\end{equation} By Theorem 5.41 in \cite[p.~225]{ln},
we have
$$|\sum_{x\in
K^*}\eta(2x+x^3)|\leq 2\sqrt{q}.$$ Therefore, if $q>9$, then the
quantity in (\ref{quantity}) is positive. The lemma now
follows.\eop
\end{pf}

Fix $e,f\in K^*$ as above. Let $H=K\times K$, $N=\{0\}\times K\leq
H$, and
$$G=\langle x,H\mid
x^4=1,(u,v)^x=(eu,fv),\forall\; (u,v)\in H\rangle,$$ where $(u,v)^x$
stands for $x^{-1}(u,v)x$. With $s_1,s_2,s_3$ as given in
Lemma~\ref{ss}, we define
\begin{equation}\label{defofRDS}
R:=R_0 +R_1x+R_2 x^2+R_3x^3\in\Bbb{Z}[G], \end{equation} where
$R_0=\{(y,y^2)\mid y\in K\}$, $R_1=\{(y,\frac{1}{s_1}y^2)\mid y\in
K\}$, $R_2=\{(y,\frac{1}{s_2}y^2)\mid y\in K\}$, and
$R_3=\{(y,\frac{1}{s_3}y^2)\mid y\in K\}$.

\begin{Theorem}\label{thm_c1}
Let $q$ be a prime power such that $q\equiv 1$ (mod 4) and $q>9$.
Then $R$ is a $(4q,q,4q,4)$ RDS in $G$ relative to $N$.
\end{Theorem}
\begin{pf}
For $(u,v)\in H$, let $\chi_{u,v}$ be the character of $H$ defined
by
$$\chi_{u,v}(u',v')=\xi_p^ {tr(uu'+vv')},\;\forall (u',v')\in H.$$ For notational convenience, we set $s_0=1$. Let $(u,v)\neq
(0,0)$. For each $i$, $0\leq i\leq 3$, we have the following facts.

Fact 1. If $v\neq 0$, then \begin{align*}
\chi_{u,v}(R_i)&=\sum_{y\in K}\xi_p^{tr(uy + \frac {v}{s_i}y^2)}\\
               &=\sum_{y\in K}\xi_p^{tr\left(\frac {v}{s_i}(y+\frac {us_i}{2v})^2 -\frac {u^2s_i}{4v}\right)}\\
               &=\eta(v)\eta(s_i)S\xi_p^{-tr(\frac{u^2s_i}{4v})}
\end{align*}\\
Fact 2. If $u\neq 0$, then $\chi_{u,0}(R_i)=\sum_{y\in K}\xi_p^{tr(uy)}=0$. \\
Fact 3. We have $\chi_{u,v}(R_i^{(-x^k)})=\overline{\chi_{e^ku,f^kv}(R_i)}$, where $R_i^{(-x^k)}=\sum_{y\in R_i}x^{-k}y^{-1}x^k$, and $k\geq 1$.\\

To prove the theorem, we will show that $RR^{(-1)}=4q+4(G-N)$, which is equivalent to the following system of group ring equations in $\Bbb{Z}[H]$:
\begin{align*}
R_0R_0^{(-1)}+R_1R_1^{(-1)}+R_2R_2^{(-1)}+R_3R_3^{(-1)}&=4q+4(H-N),\\
R_0R_1^{(-x)}+R_1R_2^{(-x)}+R_2R_3^{(-x)}+R_3R_0^{(-x)}&=4H,\\
R_0R_2^{(-x^2)}+R_2R_0^{(-x^2)}+R_1R_3^{(-x^2)}+R_3R_1^{(-x^2)}&=4H,\\
R_0R_3^{(-x^3)}+R_1R_0^{(-x^3)}+R_2R_1^{(-x^3)}+R_3R_2^{(-x^3)}&=4H.
\end{align*}
Note that the fourth equation can be obtained from the second one by
first applying $h\mapsto h^{-1}$, $\forall h\in H$, to both sides of
the second equation and then conjugating both sides of the resulting
equation by $x^3$. Therefore it suffices to show that the first
three equations hold in $\Bbb{Z}[H]$. We will do so by proving that
the left hand side and the right hand side of each of the first
three equations have the same character values for all characters of
$H$. This can be checked easily for the principal character of $H$.
Now let $\chi_{u,v}$ be an arbitrary non-principal character of $H$.
For simplicity write $\chi=\chi_{u,v}$, $\chi_1=\chi_{eu,fv}$,
$\chi_2=\chi_{e^2u,f^2v}$. Let
\begin{align*}(a,b,c,d)&=(\chi(R_0),\chi(R_1),\chi(R_2),\chi(R_3)),\\
(a',b',c',d')&=(\chi_1(R_0),\chi_1(R_1),\chi_1(R_2),\chi_1(R_3)),\\
(a'',b'',c'',d'')&=(\chi_2(R_0),\chi_2(R_1),\chi_2(R_2),\chi_2(R_3)).
\end{align*}

By Fact 3, in order to prove the theorem, it suffices to show that
\begin{align*}a\bar{a}+b\bar{b}+c\bar{c}+d\bar{d}&=4q-4\chi(N),\\
a\overline{b'}+b\overline{c'}+c\overline{d'}+d\overline{a'}&=0,\\
a\overline{c''}+c\overline{a''}+b\overline{d''}+d\overline{b''}&=0.\\
\end{align*}

If $v=0$, then $\chi$ is principal on $N$. Hence $\chi(N)=q$,
and $a=b=c=d=0$. We see that all three equations above hold in this
case.

If $v\neq 0$, then $\chi$ is non-principal on $N$. Hence
$\chi(N)=0$. Using Fact 1, we see that
$$a\bar{a}=b\bar{b}=c\bar{c}=d\bar{d}=S\overline{S}=q.$$
Therefore we have $a\bar{a}+b\bar{b}+c\bar{c}+d\bar{d}=4q-4\chi(N)$ in this case. Next we will show that
\begin{align*}
a\overline{b'}+c\overline{d'}&=0,\\
b\overline{c'}+d\overline{a'}&=0,
\end{align*}
from which it follows that $a\overline{b'}+b\overline{c'}+c\overline{d'}+d\overline{a'}=0.$ We compute $a\overline{b'}+c\overline{d'}$ as follows.
\begin{align}\label{abcd}
a\overline{b'}+c\overline{d'}&=q\eta(f)\eta(s_1)\xi_p^{-tr\left(\frac{u^2}{4v}-\frac{e^2u^2s_1}{4fv}\right)}+q\eta(f)\eta(s_2s_3)\xi_p^{-tr\left(\frac{u^2s_2}{4v}-\frac{e^2u^2s_3}{4fv}\right)}\notag\\
&=q\eta(f)\left(\eta(s_1)\xi_p^{-tr\left(\frac{u^2}{4v}-\frac{e^2u^2s_1}{4fv}\right)}+\eta(s_2s_3)\xi_p^{-tr\left(\frac{u^2s_2}{4v}-\frac{e^2u^2s_3}{4fv}\right)}\right)
\end{align}
Note that
\begin{align*}
\frac{u^2}{4v}-\frac{e^2u^2s_1}{4fv}&=\frac{u^2}{4fv}(f-e^2s_1),\\
\frac{u^2s_2}{4v}-\frac{e^2u^2s_3}{4fv}&=\frac{u^2}{4fv}(fs_2-e^2s_3).\\
\end{align*}
By the definitions of $s_1$ and $s_3$, we have
$(f-e^2s_1)=(fs_2-e^2s_3)$. Therefore,
$\frac{u^2}{4v}-\frac{e^2u^2s_1}{4fv}=\frac{u^2s_2}{4v}-\frac{e^2u^2s_3}{4fv}.$
Also, by Lemma~\ref{ss}, $\eta(s_1)=-\eta(s_2s_3)$. Combining these
two facts, we see from (\ref{abcd}) that
$a\overline{b'}+c\overline{d'}=0$. Similarly, one can show that
$b\overline{c'}+d\overline{a'}=0$. Therefore we have shown that
$a\overline{b'}+b\overline{c'}+c\overline{d'}+d\overline{a'}=0$.

To finish the proof we will show that
\begin{align*}
a\overline{c''}+b\overline{d''}&=0,\\
c\overline{a''}+d\overline{b''}&=0.
\end{align*}
We compute $a\overline{c''}+b\overline{d''}$ as follows.
\begin{align}\label{acbd}
a\overline{c''}+b\overline{d''}&=q\eta(f^2)\eta(s_2)\xi_p^{-tr\left(\frac{u^2}{4v}+\frac{u^2s_2}{4v}\right)}+q\eta(f^2)\eta(s_1s_3)\xi_p^{-tr\left(\frac{u^2s_1}{4v}+\frac{u^2s_3}{4v}\right)}\notag\\
&=q\left(\eta(s_2)\xi_p^{-tr\left(\frac{u^2}{4v}+\frac{u^2s_2}{4v}\right)}+\eta(s_1s_3)\xi_p^{-tr\left(\frac{u^2s_1}{4v}+\frac{u^2s_3}{4v}\right)}\right)
\end{align}
By the definitions of $s_1$ and $s_3$, we have $s_2+1=s_1+s_3$.
Hence
$\frac{u^2}{4v}+\frac{u^2s_2}{4v}=\frac{u^2s_1}{4v}+\frac{u^2s_3}{4v}$.
Also by Lemma~\ref{ss}, $\eta(s_2)=-\eta(s_1s_3)$. Combining these
two facts, we see from (\ref{acbd}) that
$a\overline{c''}+b\overline{d''}=0$. Similarly, we can show that
$c\overline{a''}+d\overline{b''}=0$. It follows that
$a\overline{c''}+c\overline{a''}+b\overline{d''}+d\overline{b''}=0$.
The proof is now complete.\eop
\end{pf}

{\bf Remark.} When $q=p$ is a prime, $p\equiv 1$ (mod 4), $p>9$, we
have constructed a $(4p,p,4p,4)$ RDS in groups $G_{13}'$ ($e=-f$),
$G_{14}$ ($e=1$), $G_{15}$ ($e=-1$), $G_{16}$ ($e=f$) as listed in
\cite{iia}.

\section{Non-existence of $(2p,p,2p,2)$ RDS in groups of order $2p^2$}

Throughout this section $p$ is an odd prime. We will show that there
does not exist a $(2p,p,2p,2)$ RDS in any group of order $2p^2$.

Let $G$ be a group of order $2p^2$. Then $G$ has a unique Sylow
$p$-subgroup $H$ of order $p^2$. (This is an easy consequence of
Sylow's theorems.) Hence $H$ is a normal subgroup of $G$.

We first consider the case where $H$ is cyclic. In this case, $H$
has a unique subgroup $N$ of order $p$. Hence $N$ is a normal
subgroup of $G$. Also $C_G(N)\geq H$. If $R$ is a $(2p,p,2p,2)$ RDS
in $G$ relative to $N$, then by Lemma \ref{lem_gp}, we have
$p^2|2p$, which is impossible. So from now on, we assume that $H$ is
not cyclic, say $H=\langle a,b:a^p=b^p=1,[a,b]=1\rangle$.

Let $c\in G$ be an element of order $2$. Then $G$ is a semidirect
product of $H$ and $\{1,c\}$. Since ${\rm Aut}(H)\cong
GL_2(\Bbb{F}_p)$, and every element of order 2 in $GL_2(\Bbb{F}_p)$
is conjugate to a diagonal matrix with $\pm 1$'s on the diagonal,
there are three isomorphism types of semidirect product of $H$ and
$\{1,c\}$. Below we list the three nonisomorphic groups of order
$2p^2$ with noncyclic Sylow $p$-subgroup $H$:
\begin{align*}
G_1=&\langle
a,b,c:a^p=b^p=c^2=1,[a,b]=1,a^c=a^{-1},b^c=b^{-1}\rangle;\\
G_2=&\langle a,b,c:a^p=b^p=c^2=1,[a,b]=1,a^c=a^{-1},[b,c]=1\rangle;\\
G_3=&\langle a,b,c:a^p=b^p=c^2=1,[a,b]=[a,c]=[b,c]=1\rangle.
\end{align*}

In each $G_i$, $i=1,2,3$, we consider the orbits of subgroups of
order $p$ under the action of the full automorphism group ${\rm
Aut}(G_i)$. There is only one orbit of subgroups order $p$ in $G_1$
and $G_3$, and there are three such orbits in $G_2$. We list the
orbit representatives as follows:
\begin{align*}
(1).\;\; G&=G_1, \;N=\langle a\rangle;\\
(2).\;\; G&=G_3, \;N=\langle a\rangle;\\
(3).\;\; G&=G_2, \;N=\langle a\rangle;\\
(4).\;\; G&=G_2, \;N=\langle b\rangle;\\
(5).\;\; G&=G_2, \;N=\langle ab\rangle.
\end{align*}
We remark that Case (5) is the only case where $N$ is not a normal
subgroup of $G$.

The following lemma will play an important role in our non-existence
proof.

\begin{Lemma}\label{lem1} Let $p $ be an odd prime, and let $a_0, a_1, \ldots, a_{p-1}$ be nonnegative integers such that
$\sum_{i=0}^{p-1}a_i=p$. If $A=\sum_{i=0}^{p-1}a_i\xi_p^i$ has
modulus $\sqrt{2p}$, then $p=7$, $a_s=4$, $a_{2^it+s}=1$, $0\leq
i\leq 2$, for some integers $s,t$, $0\leq s\leq 6$, $1\leq t\leq 6$,
and $a_j=0$ for the rest $j$'s.
\end{Lemma}

\begin{pf} Since $A\bar{A}=2p$, we have
$$(A)(\bar{A})=(2)(p)=(2)(1-\xi_p)^{p-1},$$
as ideals in $\Bbb{Z}[\xi_p]$. Since the ideal $(1-\xi_p)$ is fixed
by $\xi_p\mapsto \xi_p^{-1}$, we have $$(1-\xi_p)^{(p-1)/2}|(A).$$
Recall that $\Delta\bar{\Delta}=p$,
$\bar{\Delta}=(\frac{-1}{p})\Delta$, we have
$(\Delta)=(1-\xi_p)^{(p-1)/2}$. Hence $(\Delta)|(A)$, and we may
write
\begin{equation}\label{AD}
A=f(\xi_p)\Delta,
\end{equation}
where $f(\xi_p)=\sum_{i=0}^{p-1}b_i\xi_p^i$ and
$f(\xi_p)\overline{f(\xi_p)}=2$, $b_i\in \Bbb{Z}$.

Multiplying both sides of (\ref{AD}) by $\bar{\Delta}$, we have
\begin{equation}\label{AF}(\sum_{i=0}^{p-1}
a_i\xi_p^i)\left(\sum_{i=0}^{p-1}\left(\frac{-i}{p}\right)\xi_p^i\right)=p(\sum_{i=0}^{p-1}b_i\xi_p^i).\end{equation}
Comparing the coefficients of $\xi_p^k$, $k=0,1,\ldots ,(p-1)$, on
both sides of (\ref{AF}), we find that there exists some $c\in
\Bbb{Z}$ such that
$$\sum_{i}a_{k-i}\left(\frac{-i}{p}\right)=pb_k-c,\; \forall k=0,1,\ldots, (p-1).$$
Summing these equations over $k$, we get $c=\sum_{k=0}^{p-1}b_k$.
Since $(\sum_{i}b_i\xi_p^i)(\sum_{i}b_i\xi_p^{-i})=2$, we have
$$c^2=(\sum_{i=0}^{p-1}b_i)^2\equiv 2\;  \left({\rm mod}\; (1-\xi_p)\cap\Bbb{Z}\right).$$
That is, $c^2\equiv 2 $ (mod $p$). Hence $\ell:=c$ (mod $p$)$\neq
0$. Write $c=pc_1+\ell$. Note that for all $k=0,1,\ldots ,(p-1)$, on
one hand we have $|\sum_{i}a_{k-i}(\frac{-i}{p})|\leq \sum_{i\neq
0}a_{k-i}=p-a_k\leq p$, and on the other hand
$|\sum_{i}a_{k-i}\left(\frac{-i}{p}\right)|=|pb_k
-c|=|p(b_k-c_1)-\ell|$. So we must have $\delta_k:=b_k-c_1=1$ or
$0$, for all $k=0,1,\ldots ,(p-1)$. Also since
$$pc_1+\ell=\sum_{k=0}^{p-1}b_k=\sum_{k=0}^{p-1}(c_1+\delta_k),$$
we have $\sum_{k}\delta_k=\ell$. Hence exactly $\ell$ of the
$\delta_k$'s are equal to $1$. It follows that
$\sum_{k}b_k\xi_p^k=\sum_{j=1}^{\ell}\xi_p^{i_j}$. Let $S=\{i_j :
1\leq j\leq \ell\}\subset \Bbb{Z}_p$. Define
$S(x)=\sum_{j=1}^{\ell}x^{i_j}\in \Bbb{Z}[x]/(x^p-1)$. Then
$$S(x)S(x^{-1})=2+\lambda T(x),$$
where $T(x)=1+x+x^2+\cdots +x^{p-1}$, and $\lambda$ is some
nonnegative integer. It follows that $\lambda=\ell-2$ and
$\ell^2=2+\lambda p$. We then have $\lambda^2+(4-p)\lambda+2=0$.
Hence $\lambda=1$ or $2$, and $p=7$.

If $\lambda=1$, then $S$ is a $(7,3,1)$ difference set in
$\Bbb{Z}_7$. Since $2$ is a multiplier of $S$ (see
\cite[p.~323]{bjl}), we have $S=\{t+s,2t+s,4t+s\}$ for some integers
$s,t$, where $1\leq t\leq 6$. Now using $\sum_{i=0}^{6}
a_i\xi_7^{i}=\left(\sum_{i=0}^{6}(\frac{i}{7})\xi_7^i\right)(\xi_7^{t+s}+\xi_7^{2t+s}+\xi_7^{4t+s}),$
we find that there are no solutions for the $a_i$'s when
$(\frac{t}{p})=1$;  and there is a unique set of solutions:
$a_{s}=4$, $a_{2^it+s}=1$, $0\leq i\leq 2$, and $a_j=0$ for the
remaining $j$'s when $(\frac{t}{p})=-1$.

In the case where $\lambda=2$, similarly, we find that there are no
solutions for the $a_i$'s when $(\frac{t}{p})=-1$; and there is a
unique set of solutions: $a_{s}=4$, $a_{2^it+s}=1$, $0\leq i\leq 2$,
and $a_j=0$ for the remaining $j$'s, when $(\frac{t}{p})=1$.\eop
\end{pf}

We are now ready to state the main theorem in this section.

\begin{Theorem}
Let $p$ be an odd prime. Then there does not exist a $(2p,p,2p,2)$
RDS in any group of order $2p^2$.
\end{Theorem}

\begin{pf} By the analysis preceding Lemma~\ref{lem1}, we only need to
consider the five cases listed before Lemma~\ref{lem1}. We use the
same notation as in the discussion at the beginning of this section.
Suppose $R$ is a putative $(2p,p,2p,2)$ RDS in $G$ relative to $N$.
Write $R=R_1+R_2 c$, where $R_i\in \Bbb{Z}[H]$, $H=\langle
a\rangle\times\langle b\rangle\cong \Bbb{Z}_p\times \Bbb{Z}_p$. Then
$RR^{(-1)}=2p+2(G-N)$. Hence we have
\[R_1R_1^{(-1)}+R_2R_2^{(-1)}=2p+2(H-N),\quad
R_1R_2^{(-c)}+R_2R_1^{(-c)}=2H.\] Applying the principal character
of $H$ to the above equations, we find that $|R_1|=|R_2|=p$.

We now consider the five cases one by one.

\textit{Case 1.} $G=G_1$ and $N=\langle a\rangle$. In this case we
have $R_1R_2^{(-c)}=R_1R_2$. Hence
$R_1R_1^{(-1)}+R_2R_2^{(-1)}=2p+2(H-N)$ and $R_1R_2=H$. For any
$\chi\in \hat{H}$ whose restriction on $N$ is non-principal, we have
\begin{align*}
\chi(R_1)\chi(R_2)&=0,\\
\chi(R_1)\overline{\chi(R_1)}+\chi(R_2)\overline{\chi(R_2)}&=2p.
\end{align*}
Hence $|\chi(R_1)|^2=2p$ or $0$. Let $S_1=\{\chi\in \hat{H} : \chi$
is non-principal on $N$ and $|\chi(R_1)|^2=2p\}$. It is clear that
the coefficient of $1_H$ in $R_1R_1^{(-1)}$ is $|R_1|=p$. This
coefficient can also be calculated by using the inversion formula.
We therefore have
$$p=\frac {1}{p^2}\sum_{\chi\in\hat{H}}\chi(R_1R_1^{(-1)})=\frac
{1}{p^2}(p^2+2p |S_1|).$$ It follows that $|S_1|=\frac{p(p-1)}{2}$.
Now note that $Gal(\Bbb{Q}(\xi_p)/\Bbb{Q})$ acts on $\hat{H}$, and
$S_1$ is fixed (setwise) under this action. Therefore $S_1$ is
partitioned into orbits under the aforementioned action, each having
size $p-1$. So $|S_1|\equiv 0$ (mod $p-1$). But this is impossible
since
$|S_1|=\frac{p(p-1)}{2}$. We have reached the desired contradiction.\\

\textit{Case 2.} $G=G_3$ and $N=\langle a\rangle$. In this case, the
group $G$ is abelian. For any $\chi\in\hat{H}$ whose restriction to
$N$ is non-principal, we have $|\chi(R_1\pm R_2)|^2=2p$. From the
proof of Lemma~\ref{lem1}, we have $\chi(R_1+R_2)=f_1 \Delta$ and
$\chi(R_1-R_2)=f_2 \Delta$, where $f_i\in\Bbb{Z}[\xi_p]$ and
$|f_i|^2=2$, for $i=1$, 2. Since $(f_1-f_2)\Delta=2\chi(R_2)$, we
have $2|(f_1-f_2)$ in $\Bbb{Z}[\xi_p]$. Let $f_2=f_1+2x$ for some
$x\in \Bbb{Z}[\xi_p]$. Multiplying both sides of this equation by
$\bar{f_1}$, we have
$$\bar{f_1}f_2=\bar{f_1}f_1+2x\bar{f_1}=2+2x\bar{f_1}.$$
So $2|\bar{f}_1 f_2$. Let $\bar{f_1}f_2=2y$ for some $y\in
\Bbb{Z}[\xi_p]$. Multiplying both sides of the equation by $f_1$, we
obtain $f_2=f_1y$. Since both $f_1$ and $f_2$ have modulus
$\sqrt{2}$, we have $f_1=\eta f_2$ for some root of unity
$\eta\in\Bbb{Z}[\xi_p]$. Now
$2\chi(R_1)=(f_1+f_2)\Delta=f_2(1+\eta)\Delta$. Multiplying this
equation by its own complex conjugate, we find that
$2|(1+\eta)\overline{(1+\eta)}$. Recall that $\eta$ is a root of
unity in $\Bbb{Z}[\xi_p]$ and $\gcd((2), (1-\xi_p))=1$, we see that
$\eta=\pm 1$. It follows that $|\chi(R_1)|^2=0$ or $2p$. Now the
same arguments as those in the first case yield a contradiction.\\

\textit{Case 3:} $G=G_2$ and $N=\langle a\rangle$. For any $(u,v)\in
\Bbb{Z}_p^2$, we denote by $\chi_{u,v}$ the character of $H$ defined
by $\chi_{u,v}(a^{u'}b^{v'})=\xi_p^{uu'+vv'}$. Then
$\chi_{u,v}((a^ib^j)^c)=\chi_{-u,v}(a^ib^j)$. So
$\chi_{u,v}(R_i^{(-c)})=\chi_{u,-v}(R_i)$ for $i=1$, $2$. Let
$\chi\in \hat{H}$ and $\chi|_N\neq 1$. If $\chi$ is principal on
$\langle b\rangle$, then from $R_1R_2^{(-c)}+R_2R_1^{(-c)}=2H$ we
deduce that $\chi(R_1)\chi(R_2)=0$. Without loss of generality we
assume that $\chi(R_1)=0$. Then $\chi(R_2)$ has modulus $\sqrt{2p}$.
Since $R_2$ has size $p$, we have $p=7$ by Lemma \ref{lem1}. Noting
that the characters $\chi_{u,0}$ with $u\in \Bbb{Z}_p^*$ form a
single orbit of size $(p-1)$ under the action of
Gal($\Bbb{Q}(\xi_p)/\Bbb{Q}$), we have $\chi_{u,0}(R_1)=0$ for all
$u\in \Bbb{Z}_p^*$.

From $R_1R_1^{(-1)}+R_2R_2^{(-1)}=2p+2(H-N)$, we have
$R_1^{(c)}R_1^{(-c)}+R_2^{(c)}R_2^{(-c)}=2p+2(H-N)$. Now, apply a
character $\chi$ which is non-principal on $N$ to these group ring
equations, we have
\begin{align*}
&|\chi(R_1^{(c)})|^2+|\chi(R_2^{(c)})|^2=2p,\quad
|\chi(R_1)|^2+|\chi(R_2)|^2=2p,\\
&\chi(R_1)\overline{\chi(R_2^{(c)})}+\chi(R_2)\overline{\chi(R_1^{(c)})}=\chi(R_1)\chi(R_2^{(-c)})+\chi(R_2)\chi(R_1^{(-c)})=0.
\end{align*}
From the last equation, we have
$$|\chi(R_1)|^2|\chi(R_2^{(-c)}|^2=|\chi(R_2)|^2|\chi(R_1^{(-c)}|^2.$$
Substitute $|\chi(R_1)|^2$ by $2p-|\chi(R_2)|^2$, and
$|\chi(R_1^{(-c)}|^2$ by $2p-|\chi(R_2^{(-c)}|^2$ in the above
equation, we obtain
$$(2p-|\chi(R_2)|^2)|\chi(R_2^{(-c)})|^2=|\chi(R_2)|^2(2p-|\chi(R_2^{(-c)}|^2).$$
which simplifies to $|\chi(R_2)|^2=|\chi(R_2^{(-c)})|^2$. Similarly,
we can show that $|\chi(R_1)|^2=|\chi(R_1^{(-c)})|^2$. Hence
\[|\chi_{u,v}(R_i)|=|\chi_{-u,-v}(R_i)|=|\chi_{-u,v}(R_i)|=|\chi_{u,-v}(R_i)|.\]
Thus the characters of $H$ that are principal on neither $N$ nor
$\langle b\rangle$ are partitioned into subsets of size four of the
form $\{\chi_{\epsilon_1u,\epsilon_2v}:\epsilon_1,\epsilon_2=\pm
1\}$, $u,v\in \Bbb{Z}_p^*$, where $|\chi_{\epsilon_1 u,\epsilon_2
v}(R_i)|=|\chi_{u,v}(R_i)|$. Now computing the coefficient of $1_H$
in $R_1R_1^{(-1)}$ by the inversion formula, we have
$$p=\frac{1}{p^2}(p^2+\sum_{u\in \Bbb{Z}_p^*}|\chi_{u,0}(R_1)|^2+\sum_{v\in \Bbb{Z}_p^*}|\chi_{0,v}(R_1)|^2+4x)=\frac{1}{p^2}(p^2+4x)$$
for some algebraic integer $x$ . Hence $4|(p-1)$. But $p=7$: we
have reached a contradiction.\\

\textit{Case 4:} $G=G_2$ and $N=\langle b\rangle$. Let
$\chi\in\hat{H}$ and $\chi|_N\neq 1$. If $\chi$ is principal on
$\langle a\rangle$, then $\chi(R_i^{(-1)})=\chi(R_i^{(-c)})$ for
$i=1$, 2. By the same arguments as those in Case 2, we have
$|\chi(R_1)|^2=2p$ or $0$. In the former case, since $|R_1|=p$, we
have $p=7$ by Lemma~\ref{lem1}. In the latter case, we have
$|\chi(R_2)|^2=2p$. Again since $|R_2|=p$, we have $p=7$ by
Lemma~\ref{lem1}. Now the same arguments as those in the third case yield a contradiction.\\

\textit{Case 5:}  $G=G_2$ and $N=\langle ab\rangle$. Let $\chi_1$ be
the character of $H$ which maps $a$ to $1$ and $b$ to $\xi_p$. Then
$\chi_1$ is non-principal on $N$. Since $\chi_1|_{\langle
a\rangle}=1$, we have $\chi_1(R_i^{(-1)})=\chi_1(R_i^{(-c)})$. Using
the same arguments as those in Case 2, we have $|\chi_1(R_1)|^2=2p$
or $|\chi_1(R_2)|^2=2p$. Without loss of generality we assume that
$|\chi_1(R_1)|^2=2p$. Since $|R_1\cap a^iN|=1$ for all $i=0,1,\ldots
,(p-1)$, we can find a map $F_1:\Bbb{Z}_p\rightarrow\Bbb{Z}_p$ such
that $$R_1=\{a^{x+F_1(x)}b^{F_1(x)}:x\in \Bbb{Z}_p\}.$$

Let $a_i=|\{x\in \Bbb{Z}_p: F_1(x)=i\}|$. Then
$\sum_{i=0}^{p-1}a_i=p$, $a_i\geq 0$, and
$\chi_1(R_1)=\sum_{i=0}^{p-1}a_i\xi_p^i$. Since
$|\chi_1(R_1)|^2=2p$, by Lemma~\ref{lem1}, we have $p=7$,
$a_{s}=4,a_{2^it+s}=1,0\leq i\leq 2$, and $a_j=0$ for the remaining
$j$, where $s,t$ are two integers, $0\leq s\leq 6$ and $1\leq t\leq
6$. Assume that $F_1^{-1}(s)=\{i_1,i_2,i_3,i_4\}$,
$F_1^{-1}(t+s)=\{i_5\}$, $F_1^{-1}(2t+s)=\{i_6\}$,
$F_1^{-1}(4t+s)=\{i_7\}$. Now let $\chi_2$ to be the character which
maps $a$ to $\xi_p$ and $b$ to $1$.  Then
$\chi_2(R_i^{(-1)})=\chi_2(R_i)$. Combining this with
$R_1R_2^{(-c)}+R_2R_1^{(-c)}=2H$, we deduce that
$\chi_2(R_1)\chi_2(R_2)=0$. Hence $|\chi_2(R_1)|^2=0$ or $14$. That
is,
$$\chi_2(R_1)\xi_7^{-s}=(\sum_{j=1}^{4}\xi_7^{i_j})+\xi_7^{i_5+t}+\xi_7^{i_6+2t}+\xi_7^{i_7+4t}$$
has modulus $\sqrt{14}$ or $0$. We assume that $t$ is a non-square
of $\Bbb{Z}_7$. The case where $t$ is a nonzero square in
$\Bbb{Z}_7$ can be handled similarly.

We first consider the case where $|\chi_2(R_1)|^2=14$. Define
$$S(x):=(\sum_{j=1}^{4}x^{i_j})+x^{i_5+t}+x^{i_6+2t}+x^{i_7+4t}\in
\Bbb{Z}[x]/(x^7-1).$$ Then $S(x)S(x^{-1})=14+\lambda T(x)$, where
$T(x)=1+x+x^2+\cdots +x^{6}$ and $\lambda$ is a nonnegative integer.
It follows that $\lambda=\frac{1}{7}(7^2-14)=5$. Write
$S(x)=\sum_{i=0}^{6}c_i x^i$. Since the $i_j$'s are distinct, we
have $0\leq c_i\leq 4$, for all $i$. Also $\sum_{i=0}^6c_i=7$ and
$\sum_{i=0}^6c_i^2=19$. From these constrains, we find that there is
only one possibility, namely
$\{c_0,c_1,\ldots,c_6\}$=$\{4,1,1,1,0,0,0\}$. We may assume that
$c_{i_1}=4$. It follows that $i_1=i_5+t=i_6+2t=i_7+4t$. After
replacing $R_1$ by $a^{-i_1}b^{-s}R_1$ if necessary, we may assume
that $i_1=0$ and $s=0$. In order for $(\sum_{j=2}^{4}\xi_7^{i_j})+4$
to have modulus $\sqrt{14}$, we must have
$\{i_2,i_3,i_4\}=\{1,2,4\}$ or $\{3,5,6\}$ by Lemma \ref{lem1}.
Since all $i_j$'s are distinct and $t$ is assumed to be a non-square
modulo 7, we see that $\{i_2,i_3,i_4\}=\{3,5,6\}$. So $F_1$ maps all
non-squares modulo 7 to $0$, and maps each square modulo 7 to its
additive inverse. Let $\chi_3$ be the character that maps $a$ to
$\xi_7$ and $b$ to $\xi_7^u$, and $\chi_4$ be the one that maps $a$
to $\xi_7$ and $b$ to $\xi_7^{-u}$, where $u=2$ or $4$. Then it is
easy to see that $|\chi_3(R_1)|^2=7$, $|\chi_4(R_1)|^2=0$. But
similar arguments to those in Case 3 show that we must have
$|\chi_3(R_1)|=|\chi_4(R_1)|$: a contradiction.

Next we consider the case where $\chi_2(R_1)=0$. We have
$$\sum_{j=1}^{4}\xi_7^{i_j}+\xi_7^{i_5+t}+\xi_7^{i_6+2t}+\xi_7^{i_7+4t}=0.$$
Hence $\{i_1,i_2,i_3,i_4, i_5+t,i_6+2t, i_7+4t\}=\Bbb{Z}_7$. It
follows that $\{i_5,i_6,i_7\}=\{i_5+t,i_6+2t,i_7+4t\}$. Since $t\neq
0$, we have either $(i_5,i_6,i_7)=(i_5,i_5-2t,i_5+t)$ or
$(i_5,i_6,i_7)=(i_5,i_5+t,i_5+3t)$. By replacing $R_1$ with
$a^{-i_5}b^{-s}R_1$ if necessary, we may assume that $s=0$ and
$i_5=0$. When $(i_5,i_6,i_7)=(0,-2t,t)$, apply the character
$\chi_3'$ (resp. $\chi_4'$) that maps $a$ to $\xi_7$ and $b$ to
$\xi_7^u$ (resp. $\xi_7^{-u}$) to $R_1$, where $u=3$, we find that
$|\chi_3'(R_1)|^2=7$ and $|\chi_4'(R_1)|^2=0$. But again we should
have $|\chi_3'(R_1)|=|\chi_4'(R_1)|$ as before: a contradiction. The
case $(i_5,i_6,i_7)=(0,t,3t)$ is similarly ruled out: take $u=2$ (in
the definition of $\chi_3'$ and $\chi_4'$); then
$|\chi_3'(R_1)|^2=14$ and $\chi_4'(R_1)=0$, again contradicting
$|\chi_3'(R_1)|=|\chi_4'(R_1)|$.

The proof of the theorem is now complete \eop
\end{pf}

\section{Non-existence of $(4p,p,4p,4)$ RDS in abelian groups of order $4p^2$}

Throughout this section we let $G$ be an abelian group of order
$4p^2$, $p$ an odd prime. If $G$ contains a $(4p,p,4p,4)$ RDS
relative to a subgroup $N$ of order $p$, then by Lemma~\ref{lem_gp} the
Sylow $p$-subgroup of $G$ is non-cyclic. Therefore in the rest of
this section we always assume that the Sylow $p$-subgroup of $G$ is
isomorphic to $\Bbb{Z}_p\times \Bbb{Z}_p$.

In this section we will first show that if $p\neq 3$ is an odd
prime, then $G=\Bbb{Z}_2^2\times \Bbb{Z}_p^2$ does not contain a
$(4p,p,4p,4)$ RDS. We remark that $G=\Bbb{Z}_2^2\times \Bbb{Z}_3^2$
indeed contains a $(12,3,12,4)$ RDS, see \cite{djm} and \cite{mll}.

\begin{Theorem}\label{p=3}
Let $p\geq 5$ be an odd prime. Then there does not exist a
$(4p,p,4p,4)$ relative difference set in $G=\Bbb{Z}_2^2\times
\Bbb{Z}_p^2$.
\end{Theorem}

\begin{pf} We write $G=\langle\alpha_1: \alpha_1^2=1\rangle\times
\langle\alpha_2:\alpha_2^2=1\rangle\times \Bbb{Z}_p^2$ and
$H:=\Bbb{Z}_p^2<G$. Assume that $R$ is a $(4p,p,4p,4)$ RDS in $G$
relative to a subgroup $N$ of order $p$. Since the subgroups of
order $p$ of $G$ form a single orbit under the action of ${\rm
Aut}(G)$, we may choose $N$ to be $\{0\}\times \Bbb{Z}_p<H$. By the
definition of an RDS, we have
\begin{equation}\label{rdsequ}
RR^{(-1)}=4p+4(G-N)\; {\rm in}\; \Bbb{Z}[G].\end{equation} On one
hand, if $\theta\in\hat{G}$ and $\theta|_{N}=1$, then by applying
$\theta$ to both sides of (\ref{rdsequ}) we obtain that
$\theta(R)=0$. On the other hand, if $\theta\in\hat{G}$ and
$\theta|_N\neq 1$, then by applying $\theta$ to both sides of
(\ref{rdsequ}) we obtain that $\theta(R)\overline{\theta(R)}=4p$; by
the same arguments as those at the beginning of the proof of
Lemma~\ref{lem1}, we find that $\theta(R)=f_0(\xi_p)\Delta$, where
$|f_0(\xi_p)|^2=4$ and $f_0(x)\in\Bbb{Z}[x]$. Write
\begin{equation}\label{def01234}
R=R_1+R_2\alpha_1+R_3\alpha_2+R_4\alpha_1\alpha_2, \end{equation}
where $R_i\subset H$ for all $1\leq i\leq 4$. By applying the
characters of $G$ whose restrictions to $H$ are trivial to both
sides of (\ref{def01234}), we have
\begin{eqnarray}\label{trivialonN0}
|R_1|+|R_2|+|R_3|+|R_4|&=4p,\notag\\
|R_1|-|R_2|+|R_3|-|R_4|&=0,\notag\\
|R_1|-|R_2|-|R_3|+|R_4|&=0,\notag\\
|R_1|+|R_2|-|R_3|-|R_4|&=0.
\end{eqnarray}
From these equations, we find that $|R_1|=|R_2|=|R_3|=|R_4|=p$.

The characters of $H$ are of the form
$\chi_{u,v}(u',v')=\xi_p^{uu'+vv'}$, $\forall (u',v')\in H$. For any
character $\chi$ of $H$ that is non-principal on $N$, write
$(a,b,c,d)=(\chi(R_1),\chi(R_2),\chi(R_3),\chi(R_4))$. By applying
the characters of $G$ whose restrictions to $H$ equal $\chi$ to both
sides of (\ref{def01234}), we have
\begin{eqnarray}\label{f1234}
\hspace{0.1in} a+b+c+d=f_1(\xi_p)\Delta,\; \; a-b+c-d=f_2(\xi_p)\Delta,\notag\\
\hspace{0.1in} a-b-c+d=f_3(\xi_p)\Delta,\; \;
a+b-c-d=f_4(\xi_p)\Delta,
\end{eqnarray}
where $|f_i(\xi_p)|^2=4$ and $f_i(x)\in \Bbb{Z}[x]$, for
$i=1,2,3,4$. To simplify notation, we will usually write
$f_i(\xi_p)$ as $f_i$. Solving for $a,b,c,d$, we obtain,
\begin{eqnarray*}
a=\frac{1}{4}(f_1+f_2+f_3+f_4)\Delta, \; \; b=\frac{1}{4}(f_1-f_2-f_3+f_4)\Delta,\\
c=\frac{1}{4}(f_1+f_2-f_3-f_4)\Delta, \; \;
d=\frac{1}{4}(f_1-f_2+f_3-f_4)\Delta.
\end{eqnarray*}
Note that $a,b,c,d$ are all algebraic integers. We consider two
cases.

{\bf Case 1.} ${\rm ord}_p(2)$ is odd. Let $(2)=Q_1\cdots
 Q_g\bar{Q}_1\cdots\bar{Q}_g$ be the prime ideal decomposition of $(2)$ in $\Bbb{Z}[\xi_p]$.
(Note that since ${\rm ord}_p(2)$ is odd, the decomposition group of
$Q_{\ell}$ does not contain the complex conjugation.) For each $i$,
$1\leq i\leq 4$, let
$$(f_i)=Q_1^{r_{i1}}\cdots
Q_g^{r_{ig}}\bar{Q}_1^{s_{i1}}\cdots\bar{Q}_g^{s_{ig}},$$ where
$r_{i\ell}, s_{i\ell}\geq 0$. Then from $f_i\bar{f_i}=4$ we obtain
that $r_{i\ell}+s_{i\ell}=2$, $\forall \ell=1,2,\ldots g$.

We claim that $(f_i)=(f_j)$, $\forall$ $1\leq i,j\leq 4$. The proof
of the claim goes as follows. First note that by subtracting the two
equations in (\ref{f1234}) that involve $f_i$ and $f_j$, we find
that $2|(f_i-f_j)$. Hence $Q_{\ell}|(f_i-f_j)$ as well as
$\overline{Q_{\ell}}|(f_i-f_j)$ for each $\ell$. If $r_{i\ell}=0$
for some $\ell$, then $Q_{\ell}$ does not divide $f_j$ since
otherwise from $f_j\in Q_{\ell}$ and $f_i-f_j\in Q_{\ell}$ we obtain
$f_i\in Q_{\ell}$, i.e. $Q_{\ell}|(f_i)$. So we must have
$r_{j\ell}=0$. Hence $s_{i\ell}=s_{j\ell}=2$. Similarly, if
$s_{i\ell}=0$ for some $\ell$, then $s_{j\ell}=0$ and
$r_{i\ell}=r_{j\ell}=2$. If $r_{i\ell}=s_{i\ell}=1$ for some $\ell$,
then neither $r_{j\ell}$ nor $s_{j\ell}$ can be zero for otherwise
$r_{i\ell}=0$ or $2$ from the above analysis. It follows that
$r_{j\ell}=s_{j\ell}=1$. We have thus proved that $(f_i)$ and
$(f_j)$ has the same prime ideal decomposition. Hence $(f_i)=(f_j)$.
It follows that $f_i=f_1\mu_i$, where $\mu_1=1$ and $\mu_i$, $2\leq
i\leq 4$, are $2p$th roots of unity. Furthermore, if $f_i\neq \pm
f_j$ for some $i,j$, then since $(\mu_i-\mu_j)$ and $(2)$ have no
common prime ideal divisor, and $2|(f_i-f_j)$, we have
$(f_i)=(f_j)=(2)$. There are two possibilities to consider.

(i) $\mu_i=\pm 1$, $\forall i\in\{1,2,3,4\}$. In this case, noting
that $\mu_1=1$, we see that $\sum_{i=1}^4 \mu_i$ can only take one
of the values $0,4,\pm 2$. If $\sum_{i=1}^4 \mu_i=0$ or $4$, then
$(\mu_1,\mu_2,\mu_3,\mu_4)$ must be one of
\[(1,1,1,1),(1,-1,1,-1),(1,1,-1,-1),(1,-1,-1,1).\]
In each case, exactly one of $a,b,c,d$ has modulus $\sqrt{4p}$ and
the others are $0$. If $\sum_{i=1}^4 \mu_i=\pm 2$, then three of
$\mu_i$, $i=1,2,3,4$, are equal. We must have $(f_1)=(2)$ since $a$
is an algebraic integer. It follows that
$\{a,b,c,d\}=\Delta\cdot\{\eta,\eta,\eta,-\eta\}$ for some root of
unity $\eta$.

(ii) Some $\mu_i$ is not equal to $\pm 1$. In this case, we have
$(f_i)=(2)$, $\forall i\in\{1,2,3,4\}$, by the analysis immediately
preceding (i). So we write $f_i=2\omega_i$ with $\omega_i$ a root of
unity, for each $i$. It is clear that any subset of size $p-1$ of
$X:=\{1,\xi_p,\ldots,\xi_p^{p-1}\}$ forms an integral basis of
$\Bbb{Z}[\xi_p]$. So any $k$-subset of $X$ can be completed to an
integral basis of $\Bbb{Z}[\xi_p]$ when $k<p-1$. Write
$\omega_i=\epsilon_i\xi_p^{\ell_i}$ with $\ell_i\in \Bbb{Z}_p$ and
$\epsilon_i=\pm 1$ for each $i$. Then at least two of the $\ell_i
$'s are distinct, and the distinct elements among the four
$\xi_p^{\ell_i}$'s can be completed to an integral basis of
$\Bbb{Z}[\xi_p]$ as we remarked. We only consider the case where
$\ell_1\neq \ell_2$.  The remaining cases are similar. From
$a=\frac{\sum_{i=1}^4\omega_i}{2}\Delta$, we see that
$\frac{\sum_{i=1}^4\omega_i}{2}$ is an algebraic integer. Hence the
sum of coefficients of $\xi_p^{\ell_1}$  (resp. $\xi_p^{\ell_2}$) is
even in $\sum_{i=1}^4\omega_i$. Therefore we must have
$\{\ell_3,\ell_4\}=\{\ell_1,\ell_2\}$, which in turn implies that
$\{\omega_3,\omega_4\}$ is one of $\pm\{\omega_1,\omega_2\}$,
$\pm\{\omega_1,-\omega_2\}$. Case-by-case examinations show that we
must have either
$\{a,b,c,d\}=\Delta\cdot\{\eta_1,\eta_1,\eta_2,-\eta_2\}$ or
$\{a,b,c,d\}=\Delta\cdot\{\eta_1+\eta_2,\eta_1-\eta_2,0,0\}$, where
both $\eta_1$ and $\eta_2$ are roots of unity and $\eta_1\neq
\pm\eta_2$.

To summarize, we have the following three possibilities for
$(a,b,c,d)$:\\
(1A) exactly one has modulus $\sqrt{4p}$, and the others are $0$;\\
(2A) $\{a,b,c,d\}=\Delta\cdot\{\eta_1,\eta_1,\eta_2,-\eta_2\}$;\\
(2B) $\{a,b,c,d\}=\Delta\cdot\{\eta_1+\eta_2,\eta_1-\eta_2,0,0\}$,
with $\eta_1\neq \pm\eta_2$,\\
where $\eta_1,\eta_2$ are roots of unity in $\Bbb{Z}[\xi_p]$.

{\bf Case 2.} ${\rm ord}_p(2)$ is even. In this case, each prime
ideal divisor of $(2)$ in $\Bbb{Z}[\xi_p]$ is fixed by the complex
conjugation. So $f_i=2\mu_i$ for some root of unity $\mu_i$ for each
$i$. The same arguments as those in the above case work for this
case; and there are also three possibilities as listed above. In
particular, $a,b,c,d$ are multiples of $2\Delta$ in case (1A) this
time. In the following, we will consider the ${\rm ord}_p(2)$ even
case and the ${\rm ord}_p(2)$ odd case together.

First we prove that Case (2B) does not occur. Assume to the contrary
that $\chi:=\chi_{u,u'}$ with $u'\neq 0$ is a character of $H$ such
that Case (2B) occurs. Then
$\chi(R_i)=(1-\xi_p^{\ell})\xi_p^{\ell'}\epsilon\Delta$ for some
$i$, where $\ell\in \Bbb{Z}_p^*$ and $\epsilon=\pm 1$. Since $R_i$
meets each coset of $N$ in $H$ in a unique element we may write
$R_i=\{(x,f(x)):x\in \Bbb{Z}_p\}$, where $f:\Bbb{Z}_p\rightarrow
\Bbb{Z}_p$ is a function. Define $F(x):=ux+u'f(x)-\ell'$ and
$a_j:=|\{x\in\Bbb{Z}_p : F(x)=j\}|, \forall j\in \Bbb{Z}_p$. Then
\[\chi(R_i)\xi_p^{-\ell'}=\sum_{x\in\Bbb{Z}_p}\xi_p^{F(x)}=\sum_{j=0}^{p-1}
a_j\xi_p^j=(1-\xi_p^{\ell})\epsilon\Delta=\sum_{j=0}^{p-1}\bigg[\big(\frac{j}{p}\big)-\big(\frac{j-l}{p}\big)\bigg]\epsilon\xi_p^j.\]
Comparing the coefficients of $\xi^j$ on the two sides of the above
equation, we find that
$a_j-a_0=\big[(\frac{j}{p})-(\frac{j-\ell}{p})+(\frac{-\ell}{p})\big]\epsilon$.
Together with $\sum_{j=0}^{p-1}a_j=p$, we deduce that
$a_j=1+\big[(\frac{j}{p})-(\frac{j-\ell}{p})\big]\epsilon$. We now
show that there exists $j$, $1\leq j\leq p-1$, such that $a_j$ is
negative. Let $(RN)$ (resp. (NR)) be the number of pairs
$(x,x-\ell)$ in the set $1,2,\ldots,p-1$ such that $x$ (resp.
$x-\ell$) is a non-zero square modulo $p$ and $x-\ell$ (resp. $x$)
is a non-square modulo $p$. Then by elementary number theory (see,
e.g. \cite[p.~64]{ir}), we find that
\begin{eqnarray*}
(RN)=\frac{p-1}{4}+\frac {1}{2}(\delta(-\ell\in Q)-\delta(\ell\in
Q)),\\
(NR)=\frac{p-1}{4}-\frac {1}{2}(\delta(-\ell\in Q)-\delta(\ell\in
Q)),
\end{eqnarray*}
where $\delta$ is the Kronecker delta function and $Q$ is the set of
nonzero squares modulo $p$. Since $p\geq 5$, both $(RN)$ and $(NR)$
are positive. Hence there exists $j\in \Bbb{Z}_p^*$ such that
$-(\frac{j}{p})=(\frac{j-\ell}{p})=\epsilon$. It follows that
$a_j=-1<0$: a contradiction. Therefore Case (2B) can not occur.

Next we show that Case (1A) does not occur. Assume to the contrary
that $\chi:=\chi_{u,u'}$ with $u'\neq 0$ is a character of $H$ such
that Case (1A) occurs. Then $\chi(R_i)=(\sum_{j}b_j\xi_p^j)\Delta$
for some $i$, where $(\sum_{j}b_j\xi_p^j)(\sum_{j}b_j\xi_p^{-j})=4$,
$b_j\in \Bbb{Z}$. Since $R_i$ meets each coset of $N$ in $H$ in a
unique element we may write $R_i=\{(x,f(x)):x\in \Bbb{Z}_p\}$, where
$f:\Bbb{Z}_p\rightarrow \Bbb{Z}_p$ is a function. Define
$F(x):=ux+u'f(x)$ and $a_j:=|\{x\in\Bbb{Z}_p : F(x)=j\}|$. Then
$\chi(R_i)=\sum_{j} a_j\xi_p^j$. Multiplying both sides of the
following equation
$$\sum_{j}
a_j\xi_p^j=\bigg(\sum_{j}\big(\frac{j}{p}\big)\xi_p^j\bigg)\big(\sum_{j}b_j\xi_p^j\big)$$
by $\overline{\Delta}$, we get\\
\[\big(\sum_{j}
a_j\xi_p^j\big)\bigg(\sum_{j}\big(\frac{-j}{p}\big)\xi_p^j\bigg)=p(\sum_{j}b_j\xi_p^j).\]
The following arguments are similar to those in the proof of Lemma
\ref{lem1}. By comparing coefficients of $\xi_p^k$, we get
\[\sum_{j}a_{k-j}\big(\frac{-j}{p}\big)=pb_k-c,\, \forall k\in \Bbb{Z}_p\]
for some integer $c$. Summing the above equations over $k$, we get
$c=\sum_{j}b_j$. Since
$(\sum_{j}b_j\xi_p^j)(\sum_{j}b_j\xi_p^{-j})=4$, we have $c^2\equiv
4 $ mod ($(1-\xi_p)\cap \Bbb{Z}$), i.e., $c^2\equiv 4 $ (mod $p$).
Hence $c\equiv \pm 2$ (mod $p$). Write $c=pc_1+2\epsilon $ with
$\epsilon=\pm 1$. Note that
$$|p(b_k-c_1)-2\epsilon|=|pb_k-c|=|\sum_{j}a_{k-j}(\frac{-j}{p})|\leq
p-a_k\leq p.$$ So if $\epsilon=1$, then $\delta_k:=b_k-c_1=1$ or
$0$. Since $pc_1+2=\sum_{j}b_j=\sum_{j}(c_1+\delta_j) $, we have
$\sum_{j}\delta_j=2$. Hence only two of the $\delta_j$'s are equal
to $1$. It follows that $\sum_{j}b_j\xi_p^j=\xi_p^{i_1}+\xi_p^{i_2}$
with $i_1\neq i_2\in \Bbb{Z}_p$. Now $\xi_p^{i_1}+\xi_p^{i_2}$
clearly can not have modulus $2$: a
contradiction. The case where $\epsilon=-1$ is similarly ruled out.\\

So we have proved that for each character $\chi$ of $H$ that is
non-principal on $N$ only Case (2A) can possibly occur. Write
$R_i:=\{(x,h_i(x)):x\in \Bbb{Z}_p\}\subset H$ for all $i=1,2,3,4$,
where $h_i: \Bbb{Z}_p\rightarrow \Bbb{Z}_p$. For any $\chi\in
\hat{H}$ and $\chi|_N\neq 1$, we have $|\chi(R_i)|=\sqrt{p}$ for
each $i$. This implies that each $h_i$ is a $p$-ary bent function
from $\Bbb{Z}_p$ to itself. By Theorem~\ref{hou}, we have
$h_i(x)=a_i x^2+b_i x +c_i$, $a_i\neq 0$, $a_i,b_i, c_i\in
\Bbb{Z}_p$ for each $i$. For any $u\in \Bbb{Z}_p$, we write
$\chi_u:=\chi_{u,1}$, which is a character of $H$ and whose
restriction to $N$ is non-principal. Define for each $u\in\Bbb{Z}_p$
the following 4-tuple
\[(A_{1u},A_{2u},A_{3u},A_{4u})=(\chi_u(R_1),\chi_u(R_2),\chi_u(R_3),\chi_u(R_4)).\]
We have $A_{iu}=\Delta
\xi_p^{\frac{-(b_i+u)^2}{4a_i}+c_i}(\frac{a_i}{p})$ by direct
computations. Hence to meet the conditions in Case (2A), we must
have three of $(\frac{a_i}{p})$ being equal and the fourth being
distinct from them. Without loss of generality we assume that
$$\big(\frac{a_1}{p}\big)=\big(\frac{a_2}{p}\big)=\big(\frac{a_3}{p}\big)=-\big(\frac{a_4}{p}\big)$$
For each $u\in \Bbb{Z}_p$, one of the following should occur:
\begin{align*}
{\rm (i)}\;\;A_{1u}=A_{2u},A_{3u}=-A_{4u};\\
{\rm (ii)}\;\;A_{1u}=A_{3u},A_{2u}=-A_{4u};\\
{\rm (iii)}\;\;A_{2u}=A_{3u},A_{1u}=-A_{4u}.
\end{align*}
If we are in Case (i), then $a_3\neq a_4$ since
$\big(\frac{a_3}{p}\big)=-\big(\frac{a_4}{p}\big)$, and
$-\frac{(b_3+u)^2}{4a_3}+c_3=-\frac{(b_4+u)^2}{4a_4}+c_4$. The last
equation is quadratic in $u$ (the coefficient of $u^2$ is
$\frac{a_3-a_4}{4a_3a_4}\neq 0$). Therefore there are at most two
$u$'s satisfying that equation. In other words, Case (i) occurs for
at most two values of $u$. The same is true for the other two cases.
Now note that for any $u\in \Bbb{Z}_p$, one of the above three cases
must occur. It follows that $p\leq 6$. Hence $p=5$ (since $p$ is
assumed to be greater than or equal to 5). It will be convenient to
define
\begin{align*}
U_1&=\{u\in \Bbb{Z}_5\mid A_{1u}=A_{2u},A_{3u}=-A_{4u}\},\\
U_2&=\{u\in \Bbb{Z}_5\mid A_{1u}=A_{3u},A_{2u}=-A_{4u}\},\\
U_3&=\{u\in \Bbb{Z}_5\mid A_{2u}=A_{3u},A_{1u}=-A_{4u}\}.
\end{align*}
By the above analysis, we have $U_1\cup U_2\cup U_3=\Bbb{Z}_5$,
$1\leq |U_i|\leq 2$ for all $i$, and $U_i\neq U_j$ for $1\leq i\neq
j\leq 5$.

We first claim that it is impossible to have $a_1=a_2=a_3$. If
$a_1=a_2$, then $b_1\neq b_2$ since otherwise
$(a_1,b_1,c_1)=(a_2,b_2,c_2)$, which implies $U_2=U_3$, a
contradiction. Therefore, if $a_1=a_2$, then $A_{1u}=A_{2u}$ becomes
a degree one equation in $u$, which has at most one solution; hence
$|U_1|=1$. By the same reasoning we see that if $a_1=a_2=a_3$, then
$|U_1|=|U_2|=|U_3|=1$, which is clearly impossible.

Now recall that
$\big(\frac{a_1}{5}\big)=\big(\frac{a_2}{5}\big)=\big(\frac{a_3}{5}\big)$.
Since there are two non-zero squares and two nonsquares in
$\Bbb{Z}_5$, we must have two of $a_1,a_2,a_3$ being equal. Without
loss of generality assume that $a_1=a_2=-a_3$. After replacing $R$
by $R^{\sigma}g$ for some $g\in G$ and $\sigma\in Aut(G)$ which
fixes elements  of
$\langle\alpha_1\rangle\times\langle\alpha_2\rangle$, we may assume
that $h_1(x)=x^2$ (hence $a_1=1$, $b_1=c_1=0$). In the following we
study the case where $a_4=2$. The case where $a_4=-2$ can be handled
similarly.

Now that we assumed $a_1=a_2$, by the above reasoning we must have
$b_1\neq b_2$, that is $b_2\neq 0$ since $b_1$ is now assumed to be
0. We must have $|U_1|=1$, $|U_2|=|U_3|=2$, and $U_1$, $U_2$ and
$U_3$ are mutually disjoint.

Solving $A_{1u}=A_{2u}$, we see that the unique element of $U_1$ is
$u=2b_2-\frac{c_2}{2b_2}$, which must also satisfy
\begin{equation}\label{p=5i}
u^2+(-b_4-2b_3)u+2b_4^2-c_4+c_3-b_3^2=0 \end{equation} This last
equation comes from $A_{3u}=-A_{4u}$.

Any element $u\in U_2$ must satisfy \begin{align}\label{p=5ii}
2u^2+2b_3u+b_3^2-c_3=0,\\
3u^2+(-b_4+2b_2)u+2b_4^2-c_4+c_2+b_2^2=0.\end{align} Since
$|U_2|=2$, the two equations above should have two distinct common
solutions. So by comparing coefficients we have $b_4=2b_2+2b_3$ and
$2b_4^2-c_4+b_3^2-c_3+b_2^2+c_2=0$.

Now, $u=2b_2-\frac{c_2}{2b_2}\in U_1$ can not be a solution to
\eqref{p=5ii}. But adding twice of \eqref{p=5ii} to \eqref{p=5i}
gives $u=2b_2-\frac{c_2}{2b_2}$: a contradiction.

We have shown that for any character $\chi$ of $H$ that is
non-principal on $N$, none of the cases 1A, 2A, 2B can occur.
Therefore for an odd prime $p\geq 5$, a $(4p,p,4p,4)$ RDS in $G$
cannot exist. The proof is complete.\eop
\end{pf}

\begin{Theorem}\label{nop}
Let $p$ be an odd prime. Then there does not exist a $(4p,p,4p,4)$ relative difference set in $G=\Bbb{Z}_4\times \Bbb{Z}_p^2$.
\end{Theorem}
\begin{pf}
We write $G=\langle\alpha : \alpha^4=1\rangle\times \Bbb{Z}_p^2$ and
$H:=\Bbb{Z}_p^2<G$. Assume that $R$ is a $(4p,p,4p,4)$ RDS in $G$
relative to a subgroup $N$ of order $p$. Since the subgroups of
order $p$ of $G$ form a single orbit under the action of ${\rm
Aut}(G)$, we may choose $N$ to be $\{0\}\times \Bbb{Z}_p<H$. By the
definition of an RDS, we have
\begin{equation}\label{rdsequ2}
RR^{(-1)}=4p+4(G-N)\; {\rm in}\; \Bbb{Z}[G].\end{equation} On one
hand, if $\theta\in\hat{G}$ and $\theta|_{N}=1$, then by applying
$\theta$ to both sides of (\ref{rdsequ2}) we obtain that
$\theta(R)=0$. On the other hand, if $\theta\in\hat{G}$ and
$\theta|_N\neq 1$, then by applying $\theta$ to both sides of
(\ref{rdsequ2}) we obtain that $\theta(R)\overline{\theta(R)}=4p$;
by the same arguments as those at the beginning of the proof of
Lemma~\ref{lem1}, we find that $\theta(R)=f(\xi_p)\Delta$, where
$|f(\xi_p)|^2=4$ and $f(x)\in\Bbb{Z}[x]$. Write
\begin{equation}\label{def1234}
R=R_0+R_1\alpha_1+R_2\alpha^2+R_3\alpha^3, \end{equation} where
$R_j\subset H$ for $j=0,1,2$ and 3. Applying the characters of $G$
whose restrictions to $H$ are trivial to both sides of
(\ref{def1234}), we have
\begin{eqnarray}\label{trivialonN}
|R_0|+|R_1|+|R_2|+|R_3|&=4p,\notag\\
|R_0|-|R_1|+|R_2|-|R_3|&=0,\notag\\
|R_0|+i|R_1|-|R_2|-i|R_3|&=0,\notag\\
|R_0|-i|R_1|-|R_2|+i|R_3|&=0,
\end{eqnarray}
where $i^2=-1$. From these equations, we find that
$|R_0|=|R_1|=|R_2|=|R_3|=p$.

For any character $\chi\in \hat{H}$ that is non-principal on $N$,
write $(a,b,c,d)=(\chi(R_0),\chi(R_1),\chi(R_2),\chi(R_3))$. By
applying the characters of $G$ whose restrictions to $H$ are $\chi$
we obtain
\begin{align}\label{Z4abcd}
  &a+b+c+d=f_1(\xi_p)\Delta,\notag\\
   &a-b+c-d=f_2(\xi_p)\Delta,\notag\\
    &|(a-c)+(b-d)i|^2=4p,
  \end{align}
where $|f_j(\xi_p)|^2=4,\; j=1,2$, with $f_j(x)\in \Bbb{Z}[x]$. From
the first two equations in (\ref{Z4abcd}), we find that
$2(b+d)=\Delta (f_1-f_2)$. Hence $2|(f_1-f_2)$. By the same
arguments as those in the proof of Theorem~\ref{p=3}, we deduce that
$f_2=f_1\eta$ for some $2p$th root of unity $\eta\in
\Bbb{Z}[\xi_p]$. We show that $\eta$ has to be $\pm 1$. Assume to
the contrary that $\eta\neq \pm 1$. Then from $\Delta
f_1(1-\eta)=2(b+d)$ and $\Delta f_1(1+\eta)=2(a+c)$ we find that
$2|f_1$. It follows that $2|f_2$. We thus have $f_1=2\eta_1$ and
$f_2=2\eta_2$ for some roots of unity $\eta_1,\eta_2\in
\Bbb{Z}[\xi_p]$. Denote $a+c=(\eta_1+\eta_2)\Delta$ by $x$ and
$b+d=(\eta_1-\eta_2)\Delta$ by $y$. Expanding
$|(a-c)+(b-d)i|^2=|(x-2c)+(y-2d)i|^2=4p$ and noting that $1$, $i$
are linearly independent over $\Bbb{Z}[\xi_p]$, we get
\begin{align*}
x\bar{c}+\bar{x}c+y\bar{d}+\bar{y}d&=2c\bar{c}+2d\bar{d},\\
x\bar{d}+\bar{y}c-y\bar{c}-\bar{x}d-2c\bar{d}+2d\bar{c}&=p(\eta-\bar{\eta}).
\end{align*}
Here we have used the facts that $x\bar{x}+y\bar{y}=4p$ and
$x\bar{y}-\bar{x}y=2p(\eta-\bar{\eta})$. In $\Bbb{Z}[\xi_p]$, we
have $x\equiv y$ (mod 2), $\bar{x}\equiv \bar{y}$ (mod 2). So from
the above two equations we have $p(\eta-\bar{\eta})\equiv 0$ (mod
2): a contradiction. Therefore we have proved that $\eta=\pm 1$. It
follows that for an arbitrary character $\chi$ of $H$ that is
non-principal on $N$ we have
$$\chi(R_0+R_2)=0,\; |\chi(R_1+R_3)|=\sqrt{4p},$$ or
$$\chi(R_1+R_3)=0,\; |\chi(R_0+R_2)|=\sqrt{4p}.$$
We also note that for a nontrivial character $\chi$ of $H$ that is
principal on $N$ we have $\chi(R_0+R_2)=\chi(R_1+R_3)=0$ (the
argument is similar to the one we used to find $|R_j|$). By the
inversion formula, the coefficient of the identity in
$(R_0+R_2)(R_0+R_2)^{(-1)}$ is
\[\frac{1}{p^2}(4p^2+4pz)=4+\frac{4z}{p},\] where
$$z=|\{\chi\in\hat{H} : \chi|_N\neq 1, |\chi(R_0+R_2)|=\sqrt{4p}\}.$$
Hence we have $p|z$. Noting that the above set of characters is
stable under the action of Gal($\Bbb{Q}(\xi_p)/\Bbb{Q}$) on
$\hat{H}$, we see that its elements are partitioned into orbits,
each of size $p-1$. Hence $(p-1)|z$. So $z=0$ or $z=(p-1)p$. If
$z=(p-1)p$, then $\chi(R_1+R_3)=0$ for all non-principal character
$\chi$ of $H$. It follows that $R_1+R_3=\lambda H$ for some positive
integer $\lambda$. This is clearly impossible since $|R_1|=|R_3|=p$.
The case $z=0$ is similarly ruled out. The proof is complete. \eop
\end{pf}

By the analysis at the very beginning of this section, and combining
Theorem~\ref{p=3} and~\ref{nop} with the known example of a
$(12,3,12,4)$ RDS in $\Bbb{Z}_2^2\times\Bbb{Z}_3^2$ in \cite{djm} we
have

\begin{Theorem}
Let $p$ be an odd prime. An abelian group $G$ of order $4p^2$
contains a $(4p,p,4p,4)$ relative difference set if and only if
$G=\Bbb{Z}_2^2\times\Bbb{Z}_3^2$.
\end{Theorem}

\section{Conclusion}

A $(v,k,\lambda)$ difference set $D$ in a non-abelian group of order
$v$ is said to be {\it genuinely non-abelian} if none of the abelian
groups of the same order contains a difference set with these
parameters. The first genuinely non-abelian difference set was
constructed by K. Smith in \cite{ken}, and its parameters are
$(100,45,20)$.

We define a genuinely non-abelian relative difference set in the
analogous way. Combining the construction in Section 2 and the
non-existence results in Section 4, we therefore have constructed an
infinite family of genuinely non-abelian semi-regular relative
difference sets with parameters $(4p,p,4p,4)$, where $p\equiv 1$
(mod 4) is a prime and $p>9$.



\end{document}